\documentclass[10pt,twocolumn,twoside]{IEEEtran}

\IEEEoverridecommandlockouts                              % This command is only
 \pdfoutput=1
\usepackage{amsthm}
\usepackage{amsmath, amssymb}
\usepackage{amsfonts}
\usepackage{graphicx}
\usepackage{enumerate}
\usepackage{ifthen}
\usepackage{multicol}
\usepackage{float}
\usepackage{fancyhdr}
\usepackage[usenames, dvipsnames]{color}
\graphicspath{{figures/}}
\usepackage{mathtools}
\usepackage{multicol}
\usepackage{xcolor}
\usepackage{thmtools,thm-restate} % to restate theorems

\usepackage{cite}
\usepackage{svg}

 \usepackage{cancel} %To be able to cross out text
 \usepackage[normalem]{ulem}

\newtheorem{theorem}{Theorem}

\newtheorem{rem}{Remark}
\newtheorem{example}{Example}
\newtheorem{ass}{Assumption}
\newtheorem{definition}{Definition}

\newboolean{showcomments}
\setboolean{showcomments}{true}

\newcommand{\bmat}[1]{\begin{bmatrix}
#1
\end{bmatrix}}
\newcommand{\modx}[1]{\;(\mathrm{mod}\;#1)}

\newcommand{\todo}[1]{  \ifthenelse{\boolean{showcomments}}
{\textcolor{ForestGreen}{TO DO:  #1}}{}}
\newcommand{\suggest}[1]{\ifthenelse{\boolean{showcomments}}
{\textcolor{Orange}{(Suggestion: #1)}}{}}
\newcommand{\alain}[1]{\ifthenelse{\boolean{showcomments}}
{\textcolor{Blue}{(Alain says: #1)}}{}}
\newcommand{\jonas}[1]{\ifthenelse{\boolean{showcomments}}
{\textcolor{ForestGreen}{(Jonas says: #1)}}{}}
\newcommand{\kristian}[1]{\ifthenelse{\boolean{showcomments}}
{\textcolor{Blue}{(Kristian says: #1)}}{}}
\newcommand{\emma}[1]{\ifthenelse{\boolean{showcomments}}
{\textcolor{VioletRed}{(Emma says: #1)}}{}}
\newcommand{\ifneeded}[1]{\ifthenelse{\boolean{showcomments}}
{\textcolor{Gray}{#1}}{}}

\newboolean{showedit}
\setboolean{showedit}{true}
\newcommand{\edit}[1]{\ifthenelse{\boolean{showedit}}
{\textcolor{Blue}{#1}}{}}
\newcommand{\draft}[1]{\ifthenelse{\boolean{showedit}}
{\textcolor{gray}{#1}}{}}

\usepackage{array}
\newcolumntype{L}[1]{>{\raggedright\let\newline\\\arraybackslash\hspace{0pt}}m{#1}}
\newcolumntype{C}[1]{>{\centering\let\newline\\\arraybackslash\hspace{0pt}}m{#1}}
\newcolumntype{R}[1]{>{\raggedleft\let\newline\\\arraybackslash\hspace{0pt}}m{#1}}

\usepackage{physics} %to enable use of \norm

\usepackage[font = small]{caption}
\usepackage{subcaption}

\usepackage{tikz}
 \usetikzlibrary{plotmarks}
\usepackage{tikz}
\usepackage{pgfplots}
\pgfplotsset{compat=newest}
\usetikzlibrary{patterns}
\usetikzlibrary{decorations.text}
\usepgfplotslibrary{fillbetween}

\providecommand{\figref}{}
\renewcommand{\figref}[1]{Fig.~\ref{#1}}

\usepackage{lipsum}
\usepackage{mathtools}

\title{\LARGE \bf Closed-loop design for scalable performance of vehicular formations}

\author{{Jonas Hansson and Emma Tegling} 
 \thanks{The authors are with the Department of Automatic Control and the ELLIIT Strategic Research Area at Lund University, Lund, Sweden. Email: \{\tt\small{jonas.hansson, emma.tegling}\}@control.lth.se}
        \thanks{This work was partially funded by Wallenberg AI, Autonomous Systems and Software Program (WASP) funded by the Knut and Alice Wallenberg Foundation and the Swedish Research Council through Grant 2019-00691. }}

\begin{document}
\maketitle

\begin{abstract}
This paper presents a novel control design for vehicular formations, which is an alternative to the conventional second-order consensus protocol. The design is motivated by the closed-loop system, which we construct as first-order systems connected in series, and is therefore called \emph{serial consensus.} The serial consensus design will guarantee stability of the closed-loop system under the minimum requirement of the underlying communication graphs each containing a connected spanning tree --  something that is not true in general for the conventional consensus protocols. Here, we show that the serial consensus design also gives guarantees on the worst-case transient behavior of the formation, which are independent of the number of vehicles and the underlying graph structure. In particular this shows that the serial consensus design can be used to guarantee string stability of the formation, and is therefore suitable for directed formations. We show that it can be implemented through message passing or measurements to neighbors at most two hops away. The results are illustrated through numerical examples. 
\end{abstract}

\section{Introduction}\label{sec:intro}
Network systems emerge in a wide range of applications and engineered networks are, in many cases, becoming increasingly large-scale and complex. Examples include smart power grids, sensor networks, traffic and multi-robot networks, where the coordination of a multitude of interconnected subsystems or agents is a key control problem.  The prototypical coordination problem that leads to distributed consensus dynamics was studied early on by~\cite{FaxMurray2004,OlfatiSaber2004,Jadbabaie2003}, and the dynamic behaviors of this and related problems has since been the subject of much research. This has made clear that on large scales, consensus-type networks often exhibit poor dynamic behaviors, for example in terms of controllability~\cite{Pasqualetti2014}, performance and coherence~\cite{Bamieh2012Sep,SiamiMotee2015,TeglingMitra2019}, disturbance propagation~\cite{swaroop1996stringstability,seiler2004disturbancep_propagation, Besselink2018} and even instability~\cite{TEGLING2023}.  Motivated by these issues, our work proposes an alternative consensus control design with fundamentally improved scalability properties. 

We consider a classical vehicular formation control problem, in which 
each vehicle of the formation is modeled as a double integrator, whose controller relies on relative state measurements between neighboring vehicles.
In the one-dimensional case, this approach can be compactly written on the conventional second-order consensus form
\begin{equation}\label{eq:conventional_consensus} \ddot{x}(t)=u(x,t)=-L_1\dot{x}(t)-L_0x(t)+u_\mathrm{ref}(t).\end{equation}
Here, $x$ is a vector which represents the position of each vehicle, $u$ is a vector of the control inputs, $L_1$ and $L_0$ are graph Laplacians, and $u_\mathrm{ref}$ is a reference control signal. Various assumptions on the feedback structure, here captured by the two graph Laplacians, have been considered over the years. They were assumed to be proportional to each other in the early work~\cite{Ren2007Willy},  as well as in more recent analyses \cite{Patterson2014,TEGLING2023}. 
In this case, when the Laplacians capture relative and localized feedback, there are at least three problems with the design~\eqref{eq:conventional_consensus}. 
First, stability is not guaranteed for all graph Laplacians. For example, the system may be unstable if the Laplacian corresponds to a directed cycle graph, see e.g.~\cite{Studli2017,Cantos2016}. Second, in directed vehicle  strings, 
small errors may amplify throughout the formation and lead to so-called string instability --  a topic thoroughly surveyed in 
\cite{studli2017StringConcepts,FENG2019StringDefs}. Third, in the case of undirected vehicle formations, the convergence rate of the formation may scale poorly (as $O(1/N^2)$) \cite{Barooah2009}. 

A systematic study~\cite{HERMAN2017diffasym} of the feedback law \eqref{eq:conventional_consensus} noted that using different Laplacians for the position and velocity dynamics in~\eqref{eq:conventional_consensus} may drastically improve the performance of the formation. 
Specifically, symmetric position feedback, i.e. $L_0^T=L_0$, and asymmetric velocity feedback was proposed.  
However, the stability proof relies on $L_0$ being symmetric. 
In general, these systems are not straightforward to analyze in terms of the underlying topological properties, especially when considering scalability,  that is, a growth of the network. Several analytic results can, however, be derived under assumptions of spatial invariance, that is, identical agents using the same control and interacation laws \cite{Bamieh2008,Bamieh2012Sep}.

In this paper we propose a new controller for the vehicle formation: $u(x,t)=-(L_1+L_2)\dot{x}(t)-L_2L_1x(t)+u_\mathrm{ref}(t)$, where $L_1$ and $L_2$ are graph Laplacians. The controller is designed to 
give a particular closed-loop system that we call the second-order \emph{serial consensus system}. The name, as well as the reason for choosing this particular control structure, is easiest seen by considering the closed-loop system in the Laplace domain:
$$ (sI+L_2)(sI+L_1)X(s)=U_\mathrm{ref}(s).$$
Clearly, the closed-loop system has the same dynamics as two 
conventional first-order consensus systems put in a series.
As for the classical first-order consensus protocol it is true that the consensus equilibrium will be stable as long as the graphs underlying $L_1$ and $L_2$ each 
contain a connected spanning tree. This directly addresses the above mentioned problem of instability in classical second-order consensus. 
Our main results of this paper, however, concern the performance of the serial consensus. 

It turns out that the serial consensus controller can  
guarantee a strong notion of  (generalized) string stability of the formation. 
Specifically, given any combination of relative errors $e_p(t)=Lx(t)$, that are defined through a graph Laplacian~$L$, and velocity deviations $e_v(t)=\dot{x}(t)-\mathbf{1}v_\mathrm{ref}$, we are able to give bounds on the following form:
$$\sup_{t\geq 0} \left\|\bmat{e_p(t)\\e_v(t)}\right\|_\infty \leq \alpha \left\|\bmat{e_p(0)\\e_v(0)}\right\|_\infty.$$
Here $\alpha$ is a constant independent of the number of vehicles and the underlying graph structure (that is, it need not be a string, though it holds for the directed string in particular). 
The importance of this result stems from the fact that 
the bound is in terms of the $\ell^\infty$-norm. Results of this type has been suggested to be more suitable for large vehicle formations \cite{feintuch2012,studli2017StringConcepts}, especially to ensure scalability~\cite{Besselink2018}, although typically hard to derive. The result implies that our design 
addresses the earlier mentioned problems of conventional consensus. Under mild conditions the closed loop will be stable, it can be designed to achieve string stability. Since directed graphs are allowed it is possible to avoid the poor scaling of the algebraic connectivity and thus achieve a faster convergence rate \cite{Barooah2009,Hao2012stab}. 
The cost for this advantage is, in some cases, a requirement of an additional communication step, either through physical measurement or signaling. Such additional signaling has been proposed in a vehicular platooning context in e.g. \cite{Swaroop2019V2V}. 
We argue, however, that our structure gives greater benefits with a smaller communications overhead.

The remainder of the paper is organized as follows. In Section~\ref{sec:problem_setup} we introduce the problem setup and the notation used throughout the paper. Here, the serial consensus system is also defined together with some key properties. In Section~\ref{sec:main_results} we present our performance results in the form of a theorem. The results are illustrated in Section~\ref{sec:examples} through numerical examples. Finally, we conclude the paper in Section~\ref{sec:conclusions}.

\section{Problem Setup}\label{sec:problem_setup}

\subsection{Definitions and network model}
Let $\mathcal{G}=\{\mathcal{V},\mathcal{E}\}$ denote a graph of size $N=|\mathcal{V}|$ with the edge set $\mathcal{E}\subset \mathcal{V} \times \mathcal{V}$. The graph can be equivalently represented by the weighted adjacency matrix $W\in \mathbb R ^{N\times N}$ where $w_{i,j}>0 \iff (i,j)\in \mathcal{E}$. The graph is called \emph{undirected} if $W=W^T$. 
The graph contains a \emph{connected spanning tree} if for some $i\in \mathcal{V}$ there is a path from $i$ to any other vertex $j\in \mathcal{V}$.

The weighted graph Laplacian $L$ associated to the graph is defined as
\begin{equation}
    [L]_{i,j}= \left\{\begin{matrix}-w_{i,j},& ~~\mathrm{ if }\; i\neq j \\ 
                        \sum_{k\neq i} w_{i,k},& ~~\mathrm{ if }\; i=j\end{matrix}
    \right. .
    \label{eq:graph_Laplacian}
\end{equation}
Under the condition that that the graph generating the graph Laplacian contains a connected spanning tree, $L$ will have a simple and unique  eigenvalue at $0$ 
and the remaining eigenvalues will lie strictly in the right half plane (RHP). 
{We will refer to any $N\times N$ matrix that satisfies \eqref{eq:graph_Laplacian} for some set of non-negative weights $w_{i,j}$ as a graph Laplacian. }

In this work we also consider networks with a growing number of nodes. With $\mathcal{G}_N=(\mathcal{V}_N,\mathcal{E}_N)$ we denote a graph in a family~$\{\mathcal{G}_N\}$, where $N$ is the size of the growing network. 

We will denote the space of all proper, real rational, and stable transfer matrices $\mathcal{R}\mathcal{H}_\infty$ and denote the $\mathcal{H}_\infty$-norm as $\|\cdot\|_{\mathcal{H}_\infty}$ following the notation in \cite{zhou1998essentials}. By $\|\cdot\|_\infty$, we denote the standard vector norm and its corresponding induced matrix norm i.e. $\|z\|_\infty=\mathrm{max}_k|z_k|$, where $z\in \mathbb C^N$ and with $\|M\|_\infty=\sup_{\|x\|_\infty=1}\|Mx\|_\infty$, for $M\in \mathbb C^{N\times N}$.

\subsection{Vehicle formation model}
Consider a simple vehicle formation which consists of $N$ identical double integrator systems, i.e.
\begin{equation}
    \frac{\mathrm{d}^2 x_i(t)}{\mathrm{d}t^2}=u_i(x,t),\;~~ i=1,\dots,N,
    \label{eq:vehicle_formation}
\end{equation}
where $x_i(t)\in \mathbb R$.
The aim is to coordinate 
the vehicles to keep a fixed spacing and common velocity. 
This goal is related to the problem of achieving \emph{second-order consensus} as defined below.
\begin{definition}[Second-order consensus]
    The vehicle formation~\eqref{eq:vehicle_formation} is said to achieve second-order consensus if  
$$\lim_{t\rightarrow\infty}\left|\frac{\mathrm{d}x_i(t)}{\mathrm{d}t}-\frac{\mathrm{d}x_j(t)}{\mathrm{d}t}\right|=0\text{ and } \lim_{t\rightarrow\infty}|x_i(t)-x_j(t)|=0$$ for all $i,j\in \mathcal{V}$.
\end{definition}
With our control structure, the {desired, }fixed intervehicle distances can without loss of generality be set to zero for analysis purposes. This will be clarified in 
Remark~\ref{rem:pos_offsets}.

\begin{rem}\label{rem:higher_dimensions}%\jonas{Add a reference to how this can be performed.} 
In this work we only consider {a scalar state, that is, }longitudinal control. Our approach can be extended to higher spatial dimensions {(see e.g.~\cite{OhSurvey2015} for a survey of approaches)}, but we omit it here to keep notation simple.
\end{rem}

\subsection{Control structure}
In this work, we consider linear state feedback controllers of the system \eqref{eq:vehicle_formation}. Such controllers can be written as
\begin{equation}
u(t)= u_\mathrm{ref}(t)-A_1\dot{x}(t)-A_0x(t),
\label{eq:local_control}
\end{equation}
where $x=(x_1,x_2,\dots,x_N)^T\in \mathbb R^N$, $u_\mathrm{ref}(t)\in \mathbb R^N$ is a feedforward term, and $A_0,A_1\in \mathbb R^{N\times N}$ are constant feedback matrices for the position and velocity respectively. In the distributed coordination problem, the controller is further restricted to 
\begin{enumerate}[i)]
    \item only use \emph{relative} feedback;
    \item have a bounded gain;
    \item only depend on the local neighborhood of each agent.
\end{enumerate}
These restrictions are captured by considering controllers which are part of the following class.
\begin{definition}[$q$-step implementable relative feedback]\label{def:q_implementable}
    The relative state feedback $u=Ax$ is $q$-step implementable with respect to the adjacency matrix~$W$ and gain $c>0$ if $A \in\mathcal{A}^q(W,c)$,  where
    \begin{equation*}
        \mathcal{A}^q(W,c)= \begin{Bmatrix}
        A& \left|~\begin{matrix}
        \left[\sum_{k=0}^{q} W^k\right]_{i,j}=0 \implies A_{i,j}= 0,\\ 
        A \mathbf{1}=0 , ~\|A\|_\infty\leq c
        \end{matrix}
        \right.
        \end{Bmatrix}
    .\end{equation*}  
\end{definition}
Clearly, the sum of two $q$-step implementable controllers is also $q$-step implementable so if both $A_0,A_1\in \mathcal{A}^q(W,c)$ then the combined controller in \eqref{eq:local_control} will also be $q$-step implementable. To clarify the concept of $q$-step implementability, consider the following example. 
\begin{example}
    Consider a vehicle string where each agent can measure the distance to its two neighboring vehicles. This structure can be represented by the adjacency matrix $W$ such that
    $[W]_{i,j}=1 \iff |i-j|=1$ and $W_{i,j} = 0$ otherwise. 
    Then the sparsity constraint of $A\in\mathcal{A}^q(W,c)$ corresponds to the requirement that 
    $$|i-j|>q \implies [A]_{i,j}=0, $$
    i.e., that only relative measurements up to the $q$ nearest neighbors are used. One choice of a $1$-step implementable matrix ($A\in\mathcal{A}^1(W,c)$) is the graph Laplacian for an undirected path graph $A=L_\mathrm{undir-path}$ while an example of a $2$-step implementable matrix is $A=L_\mathrm{undir-path}^2$.
\end{example}

In general, if the adjacency matrix $W$ captures a physical network, then a controller $u = Ax$ with $A\in {A}^q(W,c)$ means $u_i$ only requires signals from Agent~$i$'s $q$-hop neighborhood. This is readily proven; we refer the reader to  \cite{Hansson2023ScalableConsensus}. 

A controller that has been widely applied for vehicle formations in the literature is what we will call the \emph{conventional consensus} controller. In this case both $A_0=L_0$ and $A_1=L_1$ are chosen to be graph Laplacians and this results in the controller being a $1$-step implementable relative-feedback controller. The closed-loop dynamics with this controller are
\begin{equation}
    \ddot{x}=-L_1\dot{x}-L_0x+u_\mathrm{ref}.
    \label{eq:standard_consensus}
\end{equation}
This is, however, not the only way to implement a controller satisfying the desired structure. %We propose another approach shortly.  
{We next propose our alternative approach.}

\begin{rem}\label{rem:pos_offsets} The analysis of a formation with  position offsets can be made on the translated states $\tilde{x}=x-p-tv_\mathrm{ref}\mathbf{1}$ where $p\in\mathbb{R}^N$ is a vector of desired offsets and $v_\mathrm{ref}$ a desired velocity. If the reference control signal is chosen to be $u_\mathrm{ref}=\tilde{u}_\mathrm{ref}+A_0p$ then the closed-loop dynamics in the new states becomes
$$\ddot{\tilde x}=\tilde{u}_\mathrm{ref}-A_1\dot{\tilde x}-A_0\tilde x $$
where the property $A_{0}\mathbf{1}=A_{1}\mathbf{1}=0$ %$A_{0,1} \mathbf{1}=0$ 
was used. Thus, the dynamics around any possible offset will be equivalent to the dynamics of $x$ when $p=0$ and $v_\mathrm{ref}=0$. 
\end{rem}

\subsection{A Novel Design: Serial Consensus} \label{sec:serialconsensus}
To address the stability and performance of interconnected double-integrator systems we propose a control design that achieves a desired closed loop. Due to its structure, we term it the second-order serial consensus system.
\begin{definition}[Second-order serial consensus system]\label{def:serial_consensus}
Let $L_1$ and $L_2$ be weighted and directed graph Laplacians. The second-order serial consensus system is then 
\begin{equation}
    (sI+L_2)(sI+L_1)X(s)=U_\mathrm{ref}(s).
    \label{eq:serial_consensus}
\end{equation} 
\end{definition}
{The system~\eqref{eq:vehicle_formation} achieves the serial consensus system through the control design \begin{equation}
    u(x,t)=u_\mathrm{ref}(t)-(L_2+L_1)\dot{x}(t)-L_2L_1x(t). \label{eq:serialcontroller}
\end{equation} } 

When analyzing the serial consensus controller of \eqref{eq:serial_consensus} we will make use of the following assumption on the graph structure.
\begin{ass}[Connected spanning tree]\label{ass:connected_tree}
The graphs underlying the graph Laplacians {$L_1$ and $L_2$} each contain a connected spanning tree.   
\end{ass}
A convenient state-space representation of the serial consensus system is
\begin{equation}
    \bmat{\dot{\xi}_1\\\dot{\xi}_2} \!=\!
\underbrace{\bmat{-L_1 & I\\ 0 & -L_2}}_A\! \bmat{\xi_1\\ \xi_2}+\bmat{0\\u_\mathrm{ref}}.
\label{eq:serial_xiform}
\end{equation}
This can be transformed back to $x$ and $\dot{x}$ through the linear transformations
\[\bmat{x\\ \dot{x}}=\bmat{I & 0\\ -L_1 & I} \bmat{\xi_1\\ \xi_2} ~\text{and}~ \bmat{\xi_1\\ \xi_2}=\bmat{I & 0\\ L_1 & I} \bmat{x\\ \dot{x}} \]
One benefit of considering the serial consensus can be understood through the following theorem.
\begin{restatable}{theorem}{scalstab}\label{thrm:scalable_stability}
The second-order serial consensus system~\eqref{eq:serial_consensus} under {Assumption}~\ref{ass:connected_tree} and with $U_\mathrm{ref}\in \mathcal{R}\mathcal{H}_\infty$ has the properties:
\begin{enumerate}[i)]
        \item its poles are given by the union of the eigenvalues of $-L_1$ and $-L_2$.  \label{prop:high_i}
        \item its solution achieves second-order consensus.\label{prop:high_ii}
    \end{enumerate}
\end{restatable}
The proof,  a version of which appeared in~\cite{Hansson2023ScalableConsensus}, is presented in Appendix \ref{app:thm1}. Per Theorem~\ref{thrm:scalable_stability}, the serial consensus system has a stable consensus equilibrium and, since this holds independently of the number of agents, it is also scalably stable.
For the conventional consensus a theorem like this does not exist, since using, e.g., a Laplacian corresponding to a directed cycle can result in an
unstable closed loop as noted in \cite{Studli2017,TEGLING2023}. {Serial consensus, on the other hand, can also be shown to be robustly stable, see~\cite{Hansson2023ScalableConsensus} for the robustness criteria. }

The serial consensus controller is at worst $2$-step implementable as per the following result. 
\begin{restatable}{proposition}{qimplement}\label{prop:q_implementable}
Consider the second-order serial consensus {controller~\eqref{eq:serialcontroller}.}
  %  $u(x,t)=u_\mathrm{ref}(t)-(L_2+L_1)\dot{x}(t)-L_2L_1x(t).$
    If $L_1,L_2\!\in \!\mathcal{A}^1(W,c)$ for an adjacency matrix $W$ and a constant $c$, then the controller is a %can be implemented with 
    2-step implementable relative feedback controller
    with respect to $W$ and gain $c'\!\!=\max\{2c,c^2\}$.
\end{restatable}
The proof is found in Appendix \ref{app:prop}, which is a version of what appeared in \cite{Hansson2023ScalableConsensus}. The implementation of the serial consensus will be further discussed after our main results.

\subsection{Performance criterion}
Motivated by the setting of vehicle formations we introduce the following two errors. First, let the relative position be defined by 
$$e_p(t):=d+Lx(t),$$ 
where the measurement graph $\mathcal G_N$ underlying the graph Laplacian $L$ has size $N$ and $d\in\mathbb{R}^N$ is a vector of desired offsets. Second, denote the velocity deviation as
$$e_v(t):=\dot{x}(t)-v_\mathrm{ref}\mathbf{1},$$
with $v_\mathrm{ref}\in\mathbb R$ being the desired vehicle velocity.

The error $e_p$ represents the local relative position errors. This needs to remain small to prevent vehicle collisions. Meanwhile, $e_v$ represents the deviation from the desired velocity. In a vehicle formation this error needs to be small to ensure that speed limits are respected. This should remain true also as the network grows, that is, the errors should be independent of the number of agents~$N$.

\begin{definition}[Scalable performance]
\label{def:scalable_perf}
A formation controller defined over a growing family of graphs $\{\mathcal G_N\}$ that ensures 
\begin{equation}
    \sup_{t\geq 0} \left\|\bmat{e_p(t)\\e_v(t)}\right\|_\infty \leq \alpha \left\|\bmat{e_p(0)\\e_v(0)}\right\|_\infty,
    \label{eq:scalable_performance}
\end{equation}
where $\alpha$ is fixed and independent of the number of agents $N$, is said to achieve scalable performance.
\end{definition}
Here, the choice of norm is important since the worst-case behavior gets bounded by the initial maximum deviation. We remark that we only consider the initial-value response. While the disturbance amplification scenario also requires careful analysis as discussed in~\cite{Besselink2018}, we leave it outside the scope of the present study. 

\section{Main Results}\label{sec:main_results}

Here we show that position and velocity errors can be kept small throughout the transient phase, regardless of network size. This is achieved through the use of a serial consensus controller which uses measurements based on the underlying network graph. The result is summarized in the following theorem.
\begin{theorem}    
\label{thrm:scalable_performance}
Let $e_p=Lx$ where $L$ is a graph Laplacian. If the system $\ddot{x}(t)=u(t)$ is controlled with $$u(t)=-(p_1+p_2)L\dot{x}(t)-p_1p_2L^2x(t),$$
 where $p_1,p_2>0$ and $p_1\neq p_2$, then the resulting serial consensus system achieves scalable performance with
$$\alpha =\frac{1}{|p_1-p_2|}\left(p_1+p_2+\max\{2,2p_1p_2\}\right).$$
\end{theorem}
\begin{proof}
    The serial consensus system can be rewritten as
    $$\bmat{\dot{\xi_1}\\\dot{\xi_2}}=\bmat{-p_1L & I\\0 & -p_2L }\bmat{\xi_1 \\ \xi_2}+\bmat{0\\u_\mathrm{ref}}.$$
    Here $\xi_1=x$ and $\xi_2=\dot{x}+p_1Lx$. Since $u_\mathrm{ref}=0$, the initial value problem can be solved directly. This evaluates to
    \begin{align*}
        \xi_1(t)&=e^{-p_1Lt}\xi_1(0)+e^{-p_1Lt}\int_0^te^{p_1L\tau}e^{-p_2L\tau}\mathrm{d}\tau \xi_2(0)\\
        \xi_2(t)&=e^{-p_2Lt}\xi_2(0).
    \end{align*} 
    Since $p_1L$ and $p_2L$ obviously commute, it follows that $e^{p_1L\tau}e^{-p_2L\tau}=e^{(p_1-p_2)L\tau}$. By pre-multiplying the first equation with $(p_1-p_2)L$ and using the property that $L$ commutes with $e^{-p_1Lt}$ we get the integrand \mbox{$(p_1-p_2)Le^{(p_1-p_2)L\tau}=\frac{\mathrm{d}}{\mathrm{d}t}(e^{(p_1-p_2)L\tau})$}. Finally, by applying the Fundamental Theorem of Calculus the equation can be simplified to
    $$L\xi_1\!(t)=\! e^{-p_1Lt}L\xi_1\!(0)+\frac{(e^{-p_2Lt}-e^{-p_1Lt})}{p_1-p_2}\xi_2(0).$$
    Now, inserting that $e_p(t)=L\xi_1(t)$, $\dot{x}(t)=\xi_2(t)-p_1e_p(t)$, and $e_v(t)=\dot{x}-v_\mathrm{ref}\mathbf{1}$ yields
    $$e_{p}(t)\!=\!e^{\!-p_1\!Lt}e_{p}(0)+\frac{e^{\!-p_2\!Lt}\!\!-\!e^{\!-p_1\!Lt}}{p_1-p_2}\!(e_v(t)+v_\mathrm{ref}\mathbf{1}+p_1e_{p}(0)).$$
    Next we note that the relation $e^{-Lt}\mathbf{1}=\mathbf{1}$ holds for any graph Laplacian $L$. In particular, this implies that $(e^{-p_2Lt}-e^{-p_1Lt})\mathbf{1}=0$. After some simplifications this leads to
    $$e_{p}(t)=\dfrac{p_1e^{-p_2Lt}-p_2e^{-p_1Lt}}{p_1-p_2}e_{p}(0)+\dfrac{e^{-p_2Lt}-e^{-p_1Lt}}{p_1-p_2}e_v(0)$$
    \begin{align*}        
        &e_v(t)=-p_1e_p(t)+e^{-L_2t}(e_v(0)+p_1e_p(0))\\
        &\!=\! \dfrac{p_1p_2\!\left(e^{\!-p_1Lt}\!\!-e^{\!-p_2Lt}\right)}{p_1-p_2}e_p(0)\!
        +\! \dfrac{p_1e^{\!-p_1Lt}\!\!-p_2e^{\!-p_2Lt}}{p_1-p_2}e_v(0).
    \end{align*}
    Taking the induced norm $\|\cdot\|_\infty$ and then applying the triangle inequality yields the following two bounds:
    $$\|e_p(t)\|_\infty\leq \dfrac{p_1+p_2}{|p_1-p_2|}\|e_p(0)\|_\infty +\dfrac{2}{|p_1-p_2|}\|e_v(0)\|_\infty$$
    and
    $$\|e_v(t)\|_\infty \leq \dfrac{2p_1p_2}{|p_1-p_2|}\|e_p(0)\|_\infty+\dfrac{p_1+p_2}{|p_1-p_2|}\|e_v(0)\|_\infty.$$
    Finally, we have that 
    $$\|e_p\|_\infty,\|e_v\|_\infty\leq \left\|\bmat{e_p\\e_v} \right\|_\infty=\max\{\|e_p\|_\infty,\|e_v\|_\infty\}$$
    Combining these facts yields
    $$\left\|\bmat{e_p(t)\\e_v(t)} \right\|_\infty\leq \frac{p_1+p_2+2\max \{1,p_1p_2\}}{|p_1-p_2|}\left\|\bmat{e_p(0)\\e_v(0)}\right\|_\infty,$$
    which is the definition of scalable performance with $\alpha$ as in the theorem statement.
    \end{proof}
\begin{rem}
     In the limit when $p_1$ approaches infinity and $p_2$  approaches $0$ (or vice versa) then the theoretically optimal bound of $\alpha=1$ is retrieved. For instance, let $p_1=\gamma$ and $p_2=1/\gamma$. Then, for $\gamma>1$ we get
    $$\lim_{\gamma\rightarrow \infty}\alpha(\gamma)=\lim_{\gamma\rightarrow \infty}\frac{\gamma+1/\gamma+2}{\gamma-1/\gamma}=1.$$
    In this case the dynamics are essentially reduced to a first-order consensus system, which may be desirable. This would, however, require an unbounded gain. %unlimited gain.
\end{rem}

\section{Implementation}
The serial consensus protocol has now been shown to have scalability in both stability and performance. It turns out that this also holds true for the implementation. In the vehicle formation setting we desire a controller which has a finite gain, uses few measurements, and has a decentralized implementation. 

As can be seen from~\eqref{eq:serialcontroller}, the controller must implement the graphs underlying $L_1+L_2$ and $L_2L_1$. While this requires carefully designed signaling, it is still easy to implement as a controller with finite gain and localized measurements.  

\subsection{Message passing}
\begin{figure}
    \centering
    \includegraphics[width=.95\linewidth]{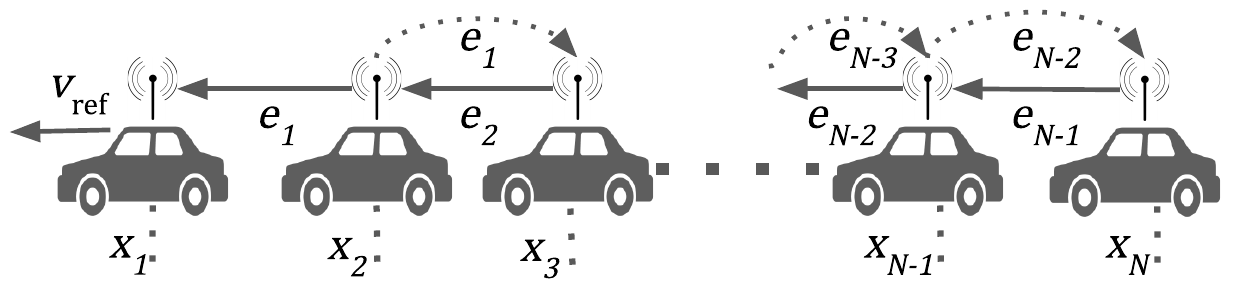}
    \caption{Vehicle platoon with directed measurements and message passing to implement serial consensus.}
    \label{fig:vehicle_platoon}
\end{figure}
One way to implement serial consensus is through message passing. The reason can be explained through the control law which, as we recall, is
$$u=-(L_2+L_1)\dot{x}-L_2L_1x+u_\mathrm{ref}.$$ 
The velocity feedback will be implementable using only relative measurements from each agent, provided that both $L_1$ and $L_2$ are localized graph Laplacians. Message passing is instead needed for the positional feedback. This can be implemented if each agent stores their relative error \mbox{$e_i = [L_1x]_i$}. 
Then if each agent can access their out-neighbor's error, then the control signal can be calculated through $L_2L_1x=L_2e$. In the case of a vehicle platoon, the relative distance to the first neighbor can be measured through the use of radar. However, the relative distance to the second neighbor requires an additional signaling layer to implement message passing. That is, if such signaling can be had with the nearest neighbor, then it is implementable in a platoon. This idea of signaling is illustrated in \figref{fig:vehicle_platoon}.
\subsection{Extended measurements}
Using a step of communication is not the only way to implement the serial consensus. Instead, $L_2L_1$ can be implemented through direct measurements. Indeed, with careful design, it is possible to choose $L_2,L_1$ so that their product $(L_2L_1)\in \mathcal{A}^1(W,c)$ and is thus implementable using only relative measurements with immediate neighbors. The following example illustrates this case.
\begin{example}
    Let $W_\mathrm{undir-path}$ correspond to the undirected path graph, i.e.
    $$(W_\mathrm{undir-path})_{i,j}=1 \iff |i-j|=1.$$
    Furthermore, let $L_\mathrm{ahead-path}$ and $L_\mathrm{behind-path}$ correspond to the look-ahead and look-behind path graphs, respectively: %, i.e.
    \begin{equation}
        \label{eq:aheadpath_def} L_\mathrm{ahead-path}=\bmat{0 & 0 \\
                                    -1 & 1 & \\
                                      & \ddots & \ddots\\
                                    & & -1 & 1} 
                                    \in \mathcal{A}^1(W_\mathrm{undir-path},2)
    \end{equation}
   and
        \begin{equation}
        \label{eq:behindpath_def} L_\mathrm{behind-path}\!\!=\!\!\bmat{1 & \!\!-1 \\
                                     & \!\!\ddots &\!\!\ddots \\
                                      &  & \!\!1 & \!\!-1\\
                                    & &  & \!\!0 & 0}\!\!\in \! \mathcal{A}^1(\!W_\mathrm{undir-path},2).
    \end{equation}
      Then, the product of these two matrices will be
    $$L_\mathrm{behind-path}L_\mathrm{ahead-path}=\bmat{1 & -1 \\
                    -1& 2 &-1 \\
                    & \ddots & \ddots& \ddots\\
                    &  & -1 & 2 & -1\\
                    & &  & & 0& 0}\!.$$
    Since the product only requires information from the neighboring states it holds that 
    $$L_\mathrm{behind-path}L_\mathrm{ahead-path}\in \mathcal{A}^1(W_\mathrm{undir-path},4),$$
    and the same holds true for $(L_\mathrm{behind-path}+L_\mathrm{ahead-path})$.  This shows that the serial consensus controller with $L_2=L_\mathrm{behind-path}$ and $L_1=L_\mathrm{ahead-path}$ would be $1$-step implementable, and thus only requires local relative feedback.
\end{example}
On the other hand, if $L_2L_1\notin\mathcal{A}^1(W,c)$, then another alternative is %might be 
to extend the local measurements. { It is then sufficient to add measurements to neighbors' neighbors. That is, given $L_1,L_2\in \mathcal{A}^1(W,c)$, then it holds that the product $(L_2L_1) \in \mathcal{A}^2(W,c')$, as guaranteed by Proposition \ref{prop:q_implementable}. The sum will clearly satisfy $(L_2+L_1)\in \mathcal{A}^1(W,2c)$. This idea of including extra measurements to improve the performance of coordination has been used in e.g.~\cite{Swaroop2019V2V} and is thus not new to the vehicular formation literature. However, in the conventional consensus, such an addition of a bounded number of neighbors does not provide improved scaling of performance in $N$, in the same way as serial consensus does. 

\begin{figure*}[!t]
    \centering
    \begin{subfigure}[]{.24\linewidth}
        \centering
        \includegraphics[width=\linewidth]{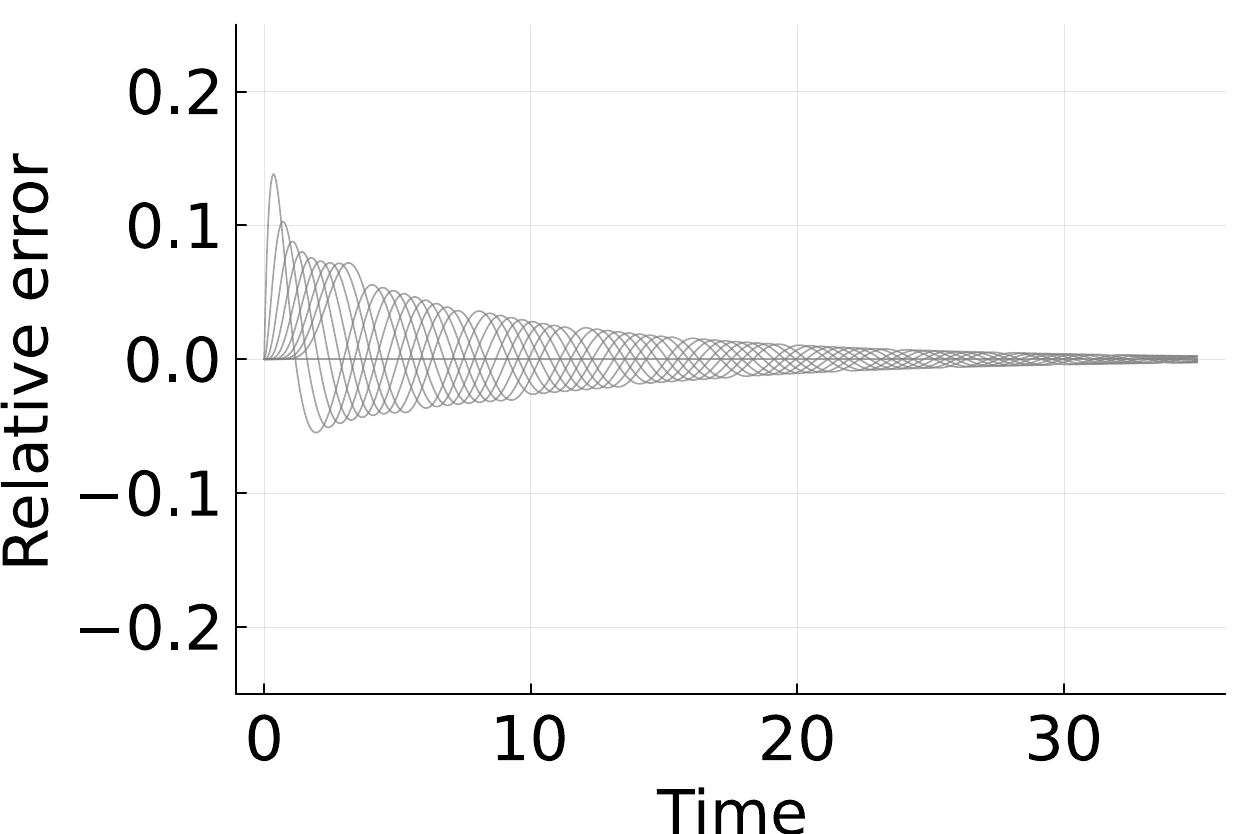}
        \caption{{Directed cycle,} conventional consensus with $N=10$.}
        \label{fig:cycle10_conv}
    \end{subfigure}   
   \hfill
   \begin{subfigure}[]{.24\linewidth}
    \centering
  \includegraphics[width=\linewidth]{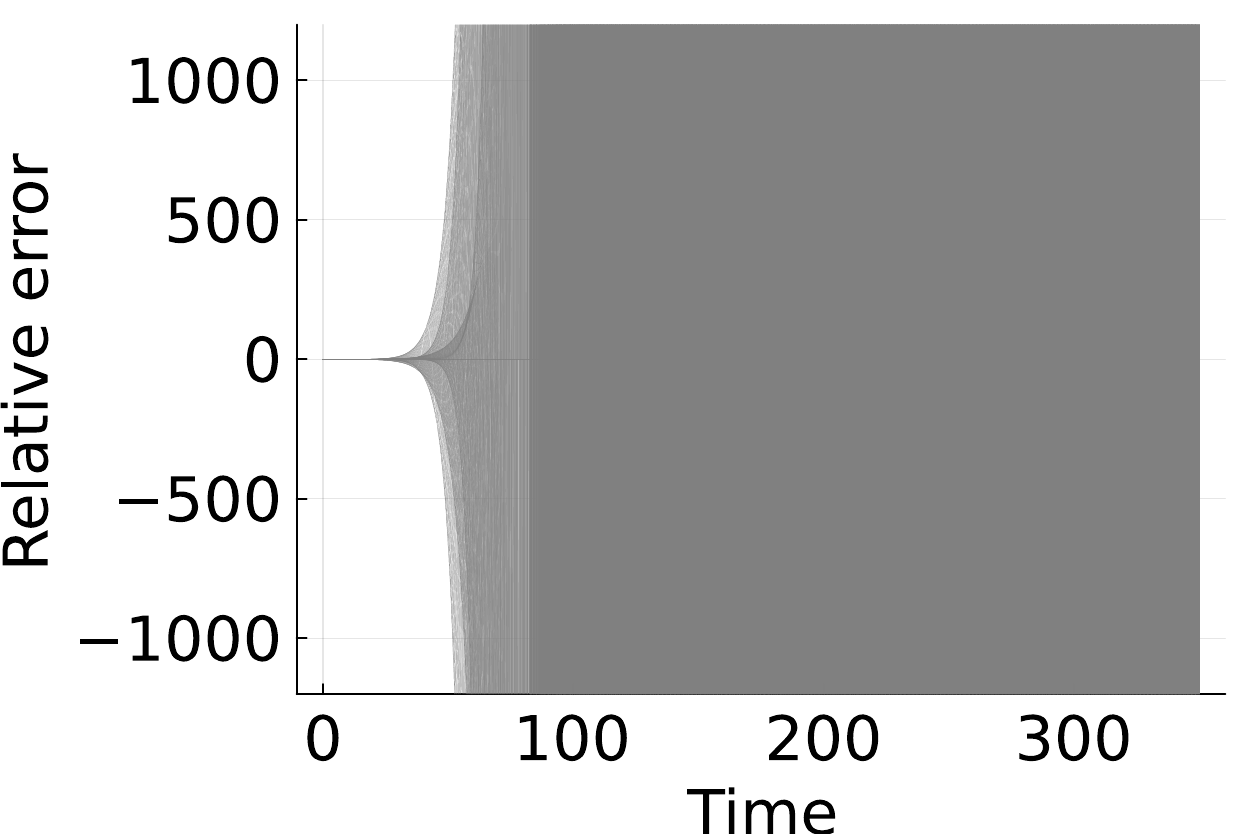}   
    \caption{{ Directed cycle, }conventional consensus with $N=100$.}
    \label{fig:cycle100_conv}
\end{subfigure}
       \hfill
   \begin{subfigure}[]{.24\linewidth}
       \centering
       \includegraphics[width=\linewidth]{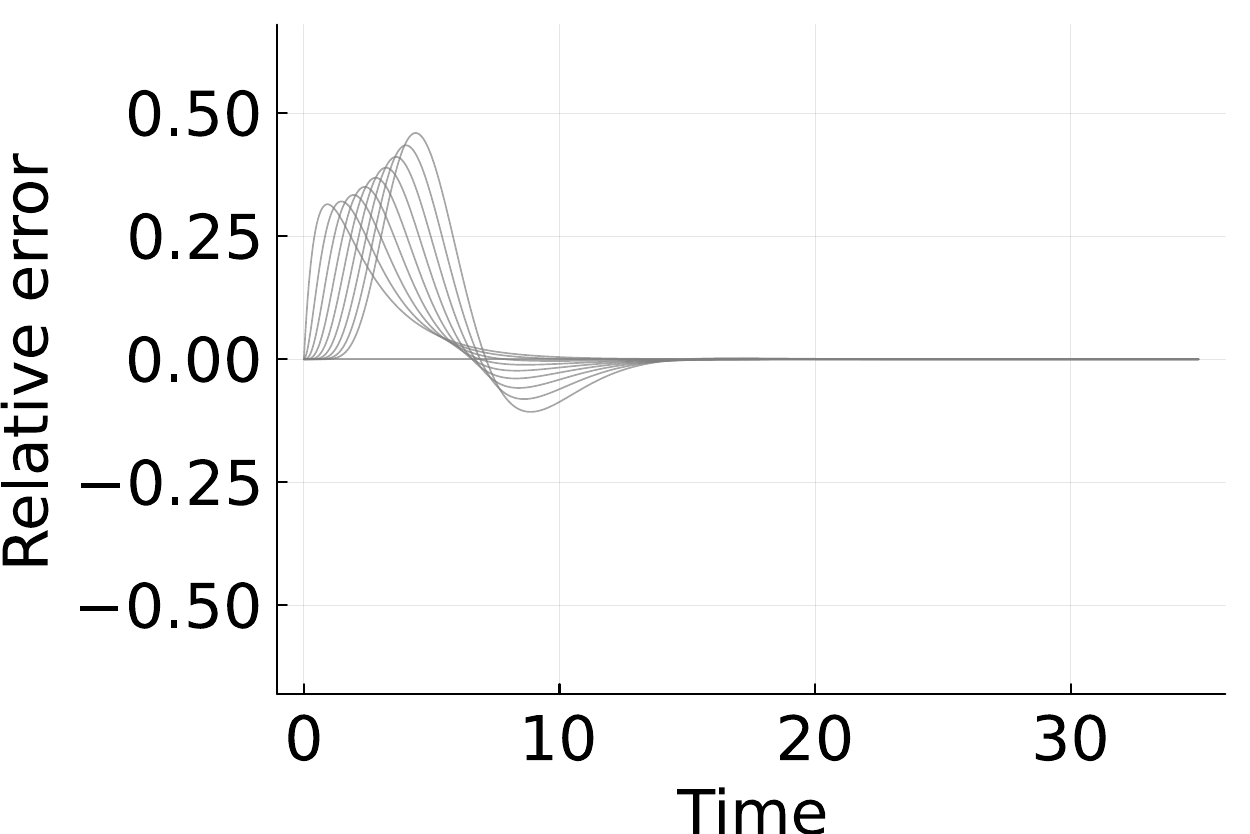}
       \caption{ Directed path, conventional consensus with $N=10$}
       \label{fig:ahead10_conv}
   \end{subfigure}
      \hfill
   \begin{subfigure}[]{.24\linewidth}
       \centering
       \includegraphics[width=\linewidth]{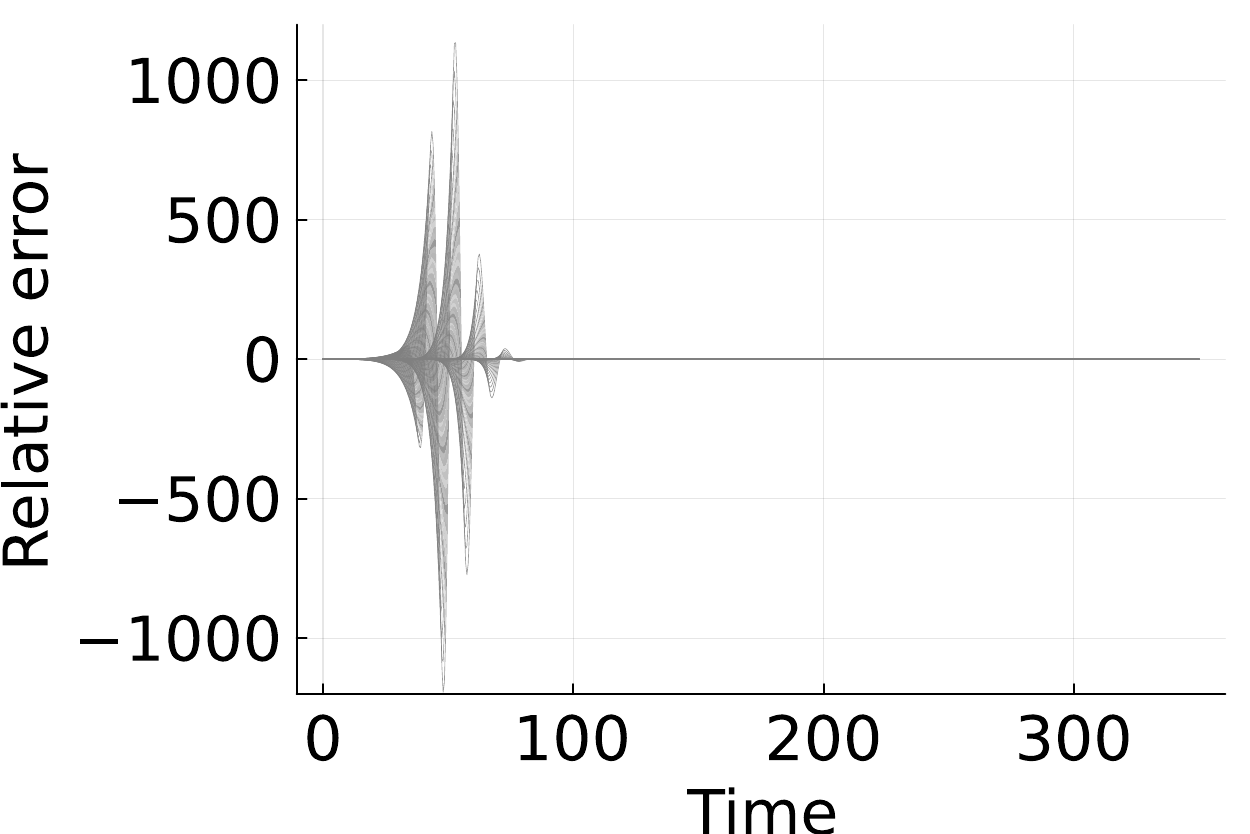}
       \caption{ Directed path, conventional consensus with $N=100$}
       \label{fig:ahead100_conv}
   \end{subfigure}
   \medskip

    \begin{subfigure}[]{.24\linewidth}
        \centering
        \includegraphics[width=\linewidth]{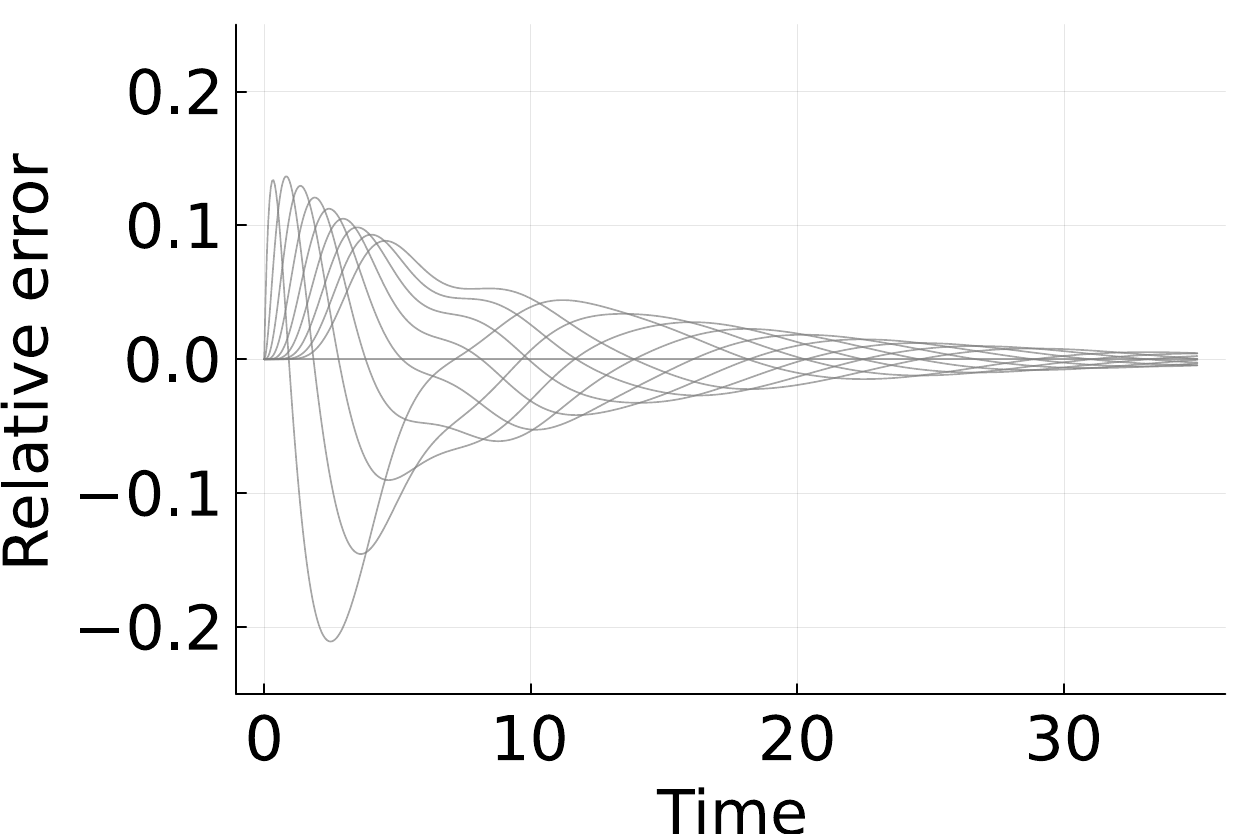}
        \caption{{Directed cycle,} serial consensus with $N=10$. }
        \label{fig:cycle10_serial}
    \end{subfigure}   
   \hfill
   \begin{subfigure}[]{.24\linewidth}
    \centering
  \includegraphics[width=\linewidth]{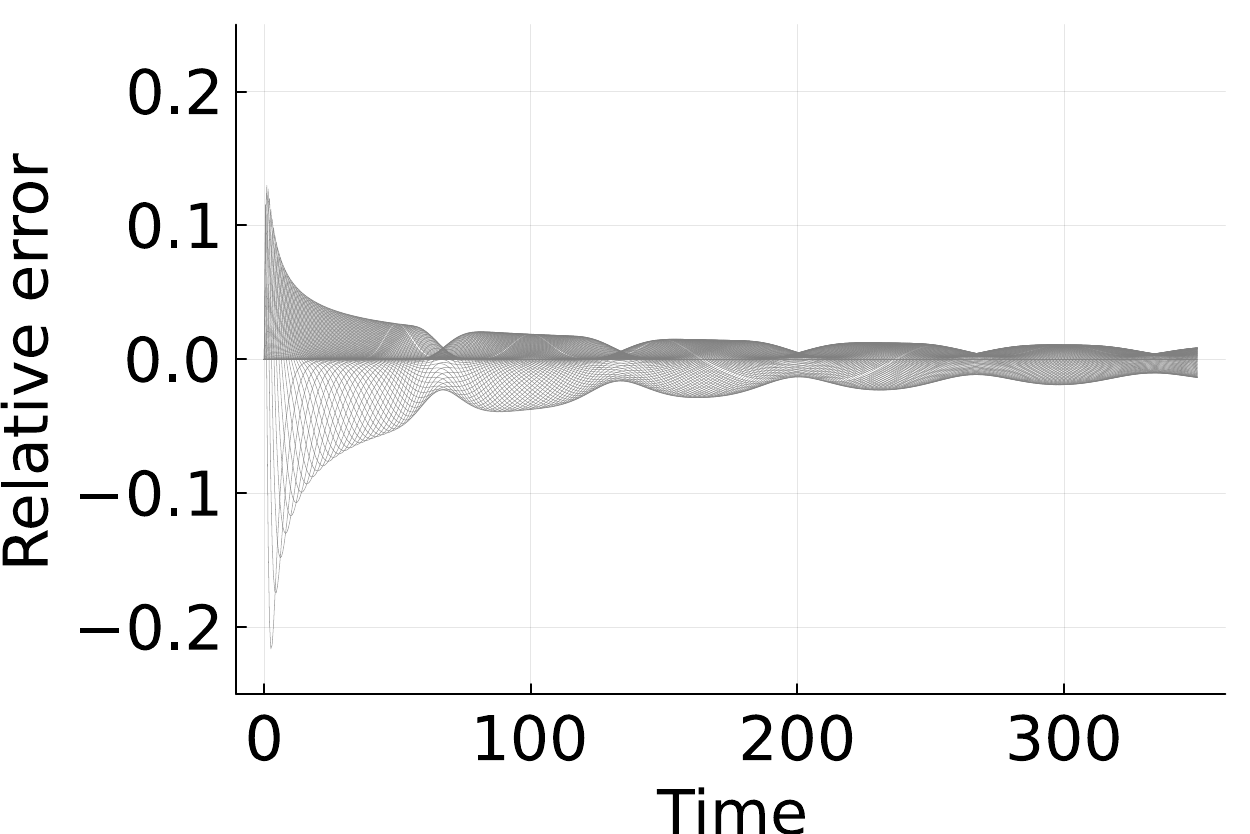}   
    \caption{{Directed cycle,} serial consensus with $N=100$.}
    \label{fig:cycle100_serial}
\end{subfigure}
       \hfill
   \begin{subfigure}[]{.24\linewidth}
       \centering
       \includegraphics[width=\linewidth]{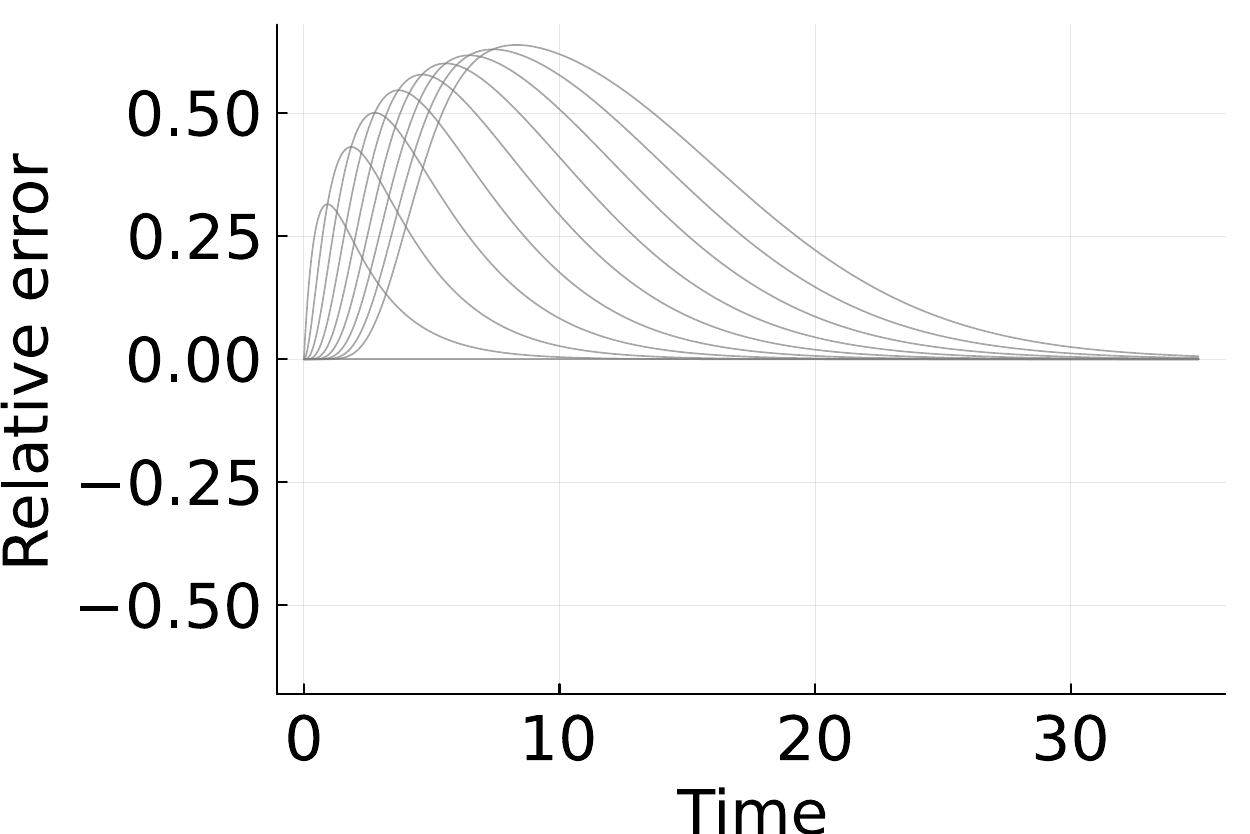}
       \caption{ {Directed path,} serial consensus with $N=10$.}
       \label{fig:ahead10_serial}
   \end{subfigure}
      \hfill
   \begin{subfigure}[]{.24\linewidth}
       \centering
       \includegraphics[width=\linewidth]{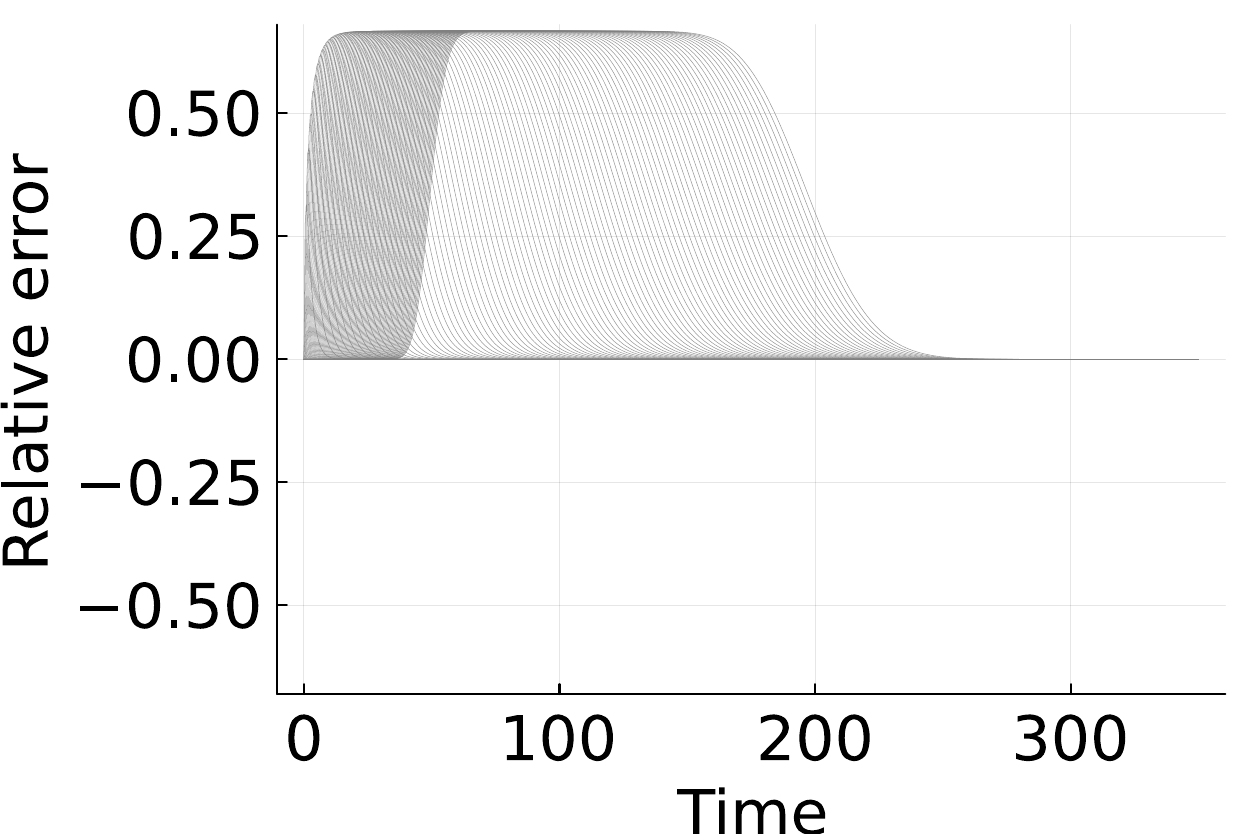}
       \caption{ {Directed path,} serial consensus with $N=100$.}
       \label{fig:ahead100_serial}
   \end{subfigure}
    \caption{Simulation of the initial value response to $x=0$, $\dot{x}_{i\neq1}(0)=0$, and $\dot{x}_1(0)=1$. For the serial consensus, $p_1=2$ and $p_2=0.5$ was used and for the conventional $r_1=2.5$ and $r_0=1$ was used. Different graph structures and number of vehicles $N$ was tested. For each plot the inter-vehicle distances $e_p(t)=L_\mathrm{ahead-path}x(t)$ are shown. The conventional consensus system can be seen to degrade with increasing number of vehicles $N$, while the serial consensus displays scalable stability and performance.}
    \label{fig:compare_serial_conv}
\end{figure*}

\section{Examples}\label{sec:examples}

In this section we will provide three examples which will illustrate our main results and how  serial consensus compares to the conventional consensus protocol. 

\subsection{Scalable stability}
\begin{example}

Consider the uni-directional circular graph structure with the graph Laplacian defined as
$$(L_\mathrm{ahead-cycle})_{i,j}=\begin{cases}
 1& \text{ if } i=j \\ 
-1& \text{ if } i=j+1 \modx{N}
\end{cases}.$$
The conventional consensus protocol is then
$$\ddot{x}=-r_1L_\mathrm{ahead-cycle}\dot x -r_0L_\mathrm{ahead-cycle}x+ u_\mathrm{ref}$$
This system is known to be troublesome, since unless $r_0$ and $r_1$ are chosen to depend on the number of vehicles $N$, the closed-loop system will be unstable for large $N$. Algebraically, this is a consequence of the smallest in magnitude eigenvalues of $L_\mathrm{ahead-cycle}$ approaching the origin at an angle to the real axis as $N$ grows \cite{TEGLING2023}. 

On the other hand, the serial consensus protocol 
$$ \ddot{x}= -(p_1+p_2)L_{\mathrm{ahead-cycle}}\dot x -p_1p_2L^2_{\mathrm{ahead-cycle}}x+ u_\mathrm{ref}$$
is stable for any $N$ as long as $p_1$ and $p_2$ are chosen to be positive, which follows from Theorem \ref{thrm:scalable_stability}. Provided $p_1\neq p_2$, Theorem \ref{thrm:scalable_performance} also asserts that it has scalable performance with respect to the graph Laplacian $L_\mathrm{ahead-cycle}$.

A comparison of the transient behavior for the two formations is shown in \figref{fig:compare_serial_conv}. The figure shows how the formation that is controlled through the serial consensus protocol is stable and has similar behavior independent of the number of vehicles. Meanwhile, the one controlled with conventional consensus has similar performance for small $N$ but eventually loses stability for large $N$.

\end{example}

\subsection{Scalable performance}
\begin{example}
Here we will illustrate the significant difference in performance between conventional and serial consensus in the case of a directed vehicle string. 
 For this purpose, consider the directed path topology whose graph Laplacian is given by~\eqref{eq:aheadpath_def}.
In this case it is easy to verify that the { conventional }%regular 
consensus protocol 
$$\ddot{x}=-r_1L_\mathrm{ahead-path}\dot x -r_0L_\mathrm{ahead-path}x+ u_\mathrm{ref}$$
will stabilize the vehicle formation for any choice of positive~$r_1$ and $r_0$. However, the transient behavior will scale poorly independent of the choice of $r_1$ and $r_0$ as is illustrated in \figref{fig:compare_serial_conv}. { The formation under this control is well known to lack string stability \cite{seiler2004disturbancep_propagation}, that is, disturbances propagate and grow along the string. }On the other hand, the serial consensus protocol with 
$$\ddot{x}=-(p_1+p_2)L_\mathrm{ahead-path}\dot x -p_1p_2L_\mathrm{ahead-path}^2x+ u_\mathrm{ref}$$
will have scalable performance with respect to the {position error $e_p = d+ L_\mathrm{ahead-path}x$ and velocity error $e_v$}
as long as $p_1\neq p_2$. This is illustrated in \figref{fig:ahead10_serial} and \figref{fig:ahead100_serial}, where the same initial conditions and parameters are used for both %the regular and serial 
consensus protocols. From the figures it is clear that both formations are stable. However, the conventional consensus protocol has a much larger transient than the serial consensus protocol, which gets worse as the number of vehicles increases.
\end{example}

\subsection{Different graph Laplacians}
\begin{example}
    The conventional consensus design has been shown to have {acceptable} performance {in a vehicle string }when different Laplacians are used {in the position and velocity feedback }\cite{HERMAN2017diffasym}. In particular, the use of the directed {path} graph Laplacian for the velocity term {(look-ahead)} and an undirected {path} graph Laplacian for the positional term {(look-ahead and look-behind)}, defined as 
    $$(L_\mathrm{undir-path})_{i,j}\!=\!\begin{cases}
        2& \text{ if } i=j \text{ and } 2\leq i\leq N-1 \\ 
        1& \text{ if } i=j \text{ for } i= 1, N \\ 
       -1& \text{ if } |i-j|=1.
       \end{cases}$$
    The resulting closed-loop system is then
    \begin{equation}\label{eq:undir_conv}\ddot{x}=\!\!-r_1L_\mathrm{ahead-path}\dot x -r_0L_\mathrm{undir-path}x+\! u_\mathrm{ref}(t).\end{equation}
    The step responses for $N=10,100$ can be seen in \figref{fig:undir10_conv}~and~\ref{fig:undir100_conv}. 
    
    This can be compared to the serial consensus utilizing bidirectional information. For instance, if the forward-looking graph Laplacian $L_\mathrm{ahead-path}$ is used together with the corresponding backward-looking graph Laplacian 
    $$(L_\mathrm{behind-path})_{i,j}=\begin{cases}
        1& \text{ if } i=j \text{ and } i\leq N-1 \\ 
       -1& \text{ if } i=j-1
       \end{cases}.$$
    The resulting closed-loop system is then
    \begin{multline}\label{eq:undir_serial}\ddot{x}(t)=-(p_1L_\mathrm{ahead-path}+p_2L_\mathrm{behind-path})\dot x(t)\\ -p_1p_2L_\mathrm{ahead-path}L_\mathrm{behind-path}x(t)+ u_\mathrm{ref}(t).
    \end{multline}
    The step responses can be seen in \figref{fig:bidirec_serial}. From the figure we can observe that the serial consensus and conventional consensus can have similar transient performance. We find it is, however, easier to predict that the serial consensus will perform well, than various versions of the conventional protocol. Indeed, the protocol proposed by~\cite{HERMAN2017diffasym} is very similar to the serial consensus controller in Example~4. 
\end{example}

\begin{figure}[]
    \centering
    \begin{subfigure}[]{.48\linewidth}
        \centering
        \includegraphics[width=\linewidth]{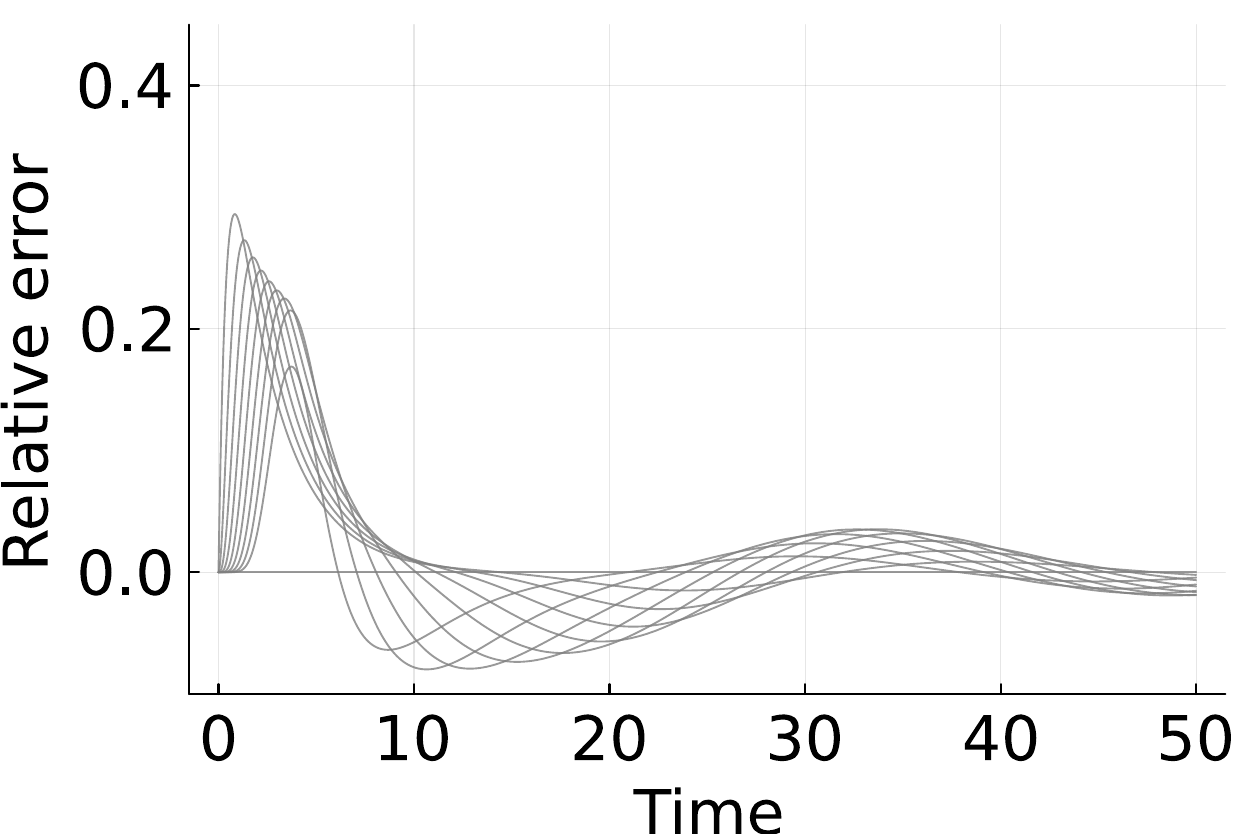}
        \caption{Conventional { consensus with }symmetric position, asymmetric velocity { feedback} and $N=10$.}
        \label{fig:undir10_conv}
    \end{subfigure}   
   \hfill
   \begin{subfigure}[]{.48\linewidth}
    \centering
  \includegraphics[width=\linewidth]{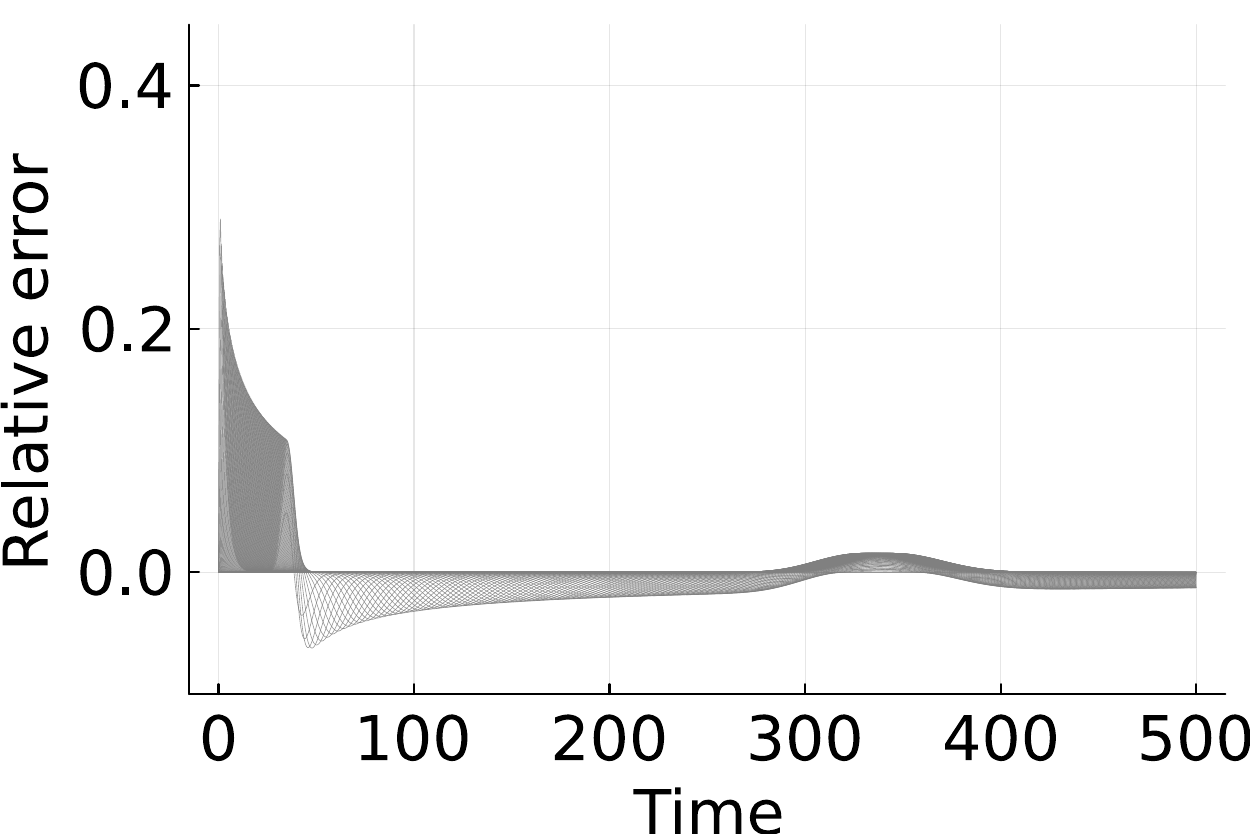}
    \caption{Conventional { consensus with }symmetric position, asymmetric velocity { feedback} and $N=100$.}
    \label{fig:undir100_conv}
\end{subfigure}
   \medskip

    \begin{subfigure}[]{.48\linewidth}
        \centering
        \includegraphics[width=\linewidth]{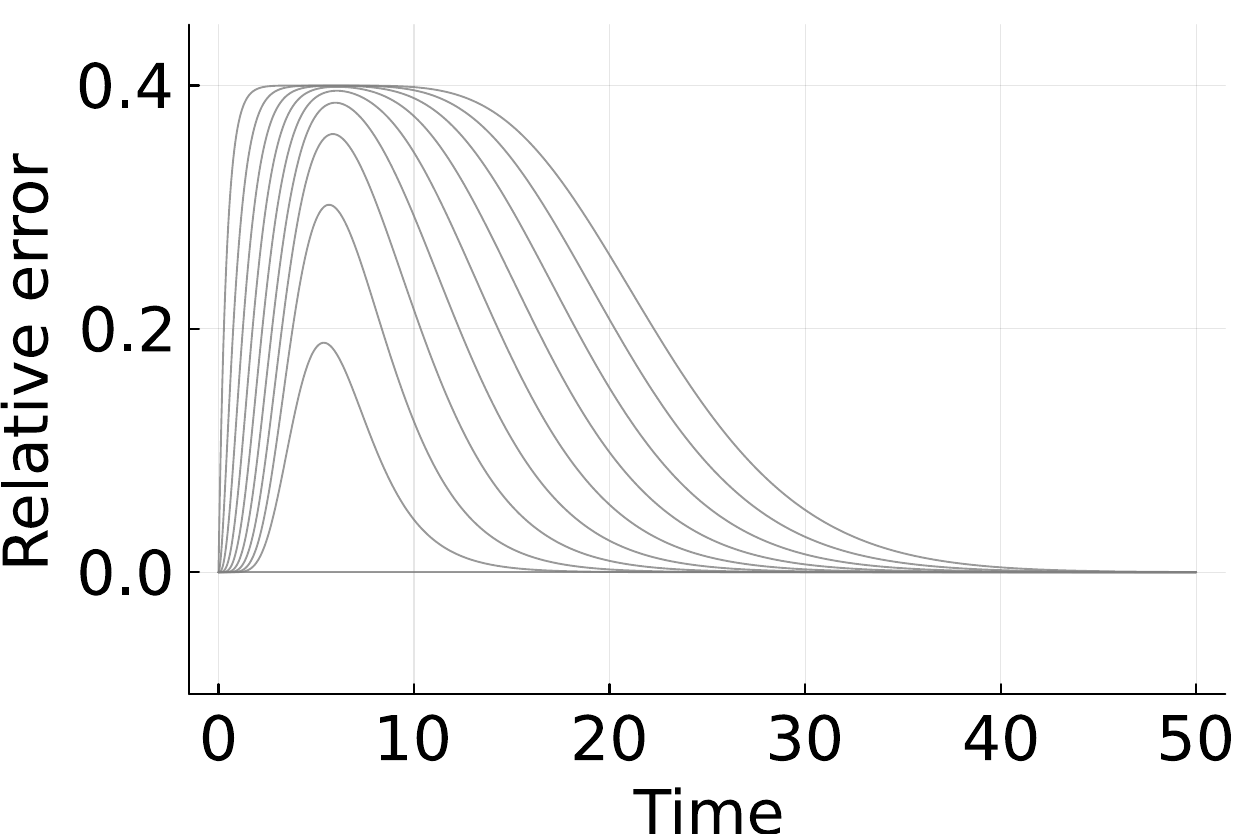}
        \caption{Serial { consensus with bidirectional} feedback and $N=10$. }
        \label{fig:undir10_serial}
    \end{subfigure}   
   \hfill
   \begin{subfigure}[]{.48\linewidth}
    \centering
  \includegraphics[width=\linewidth]{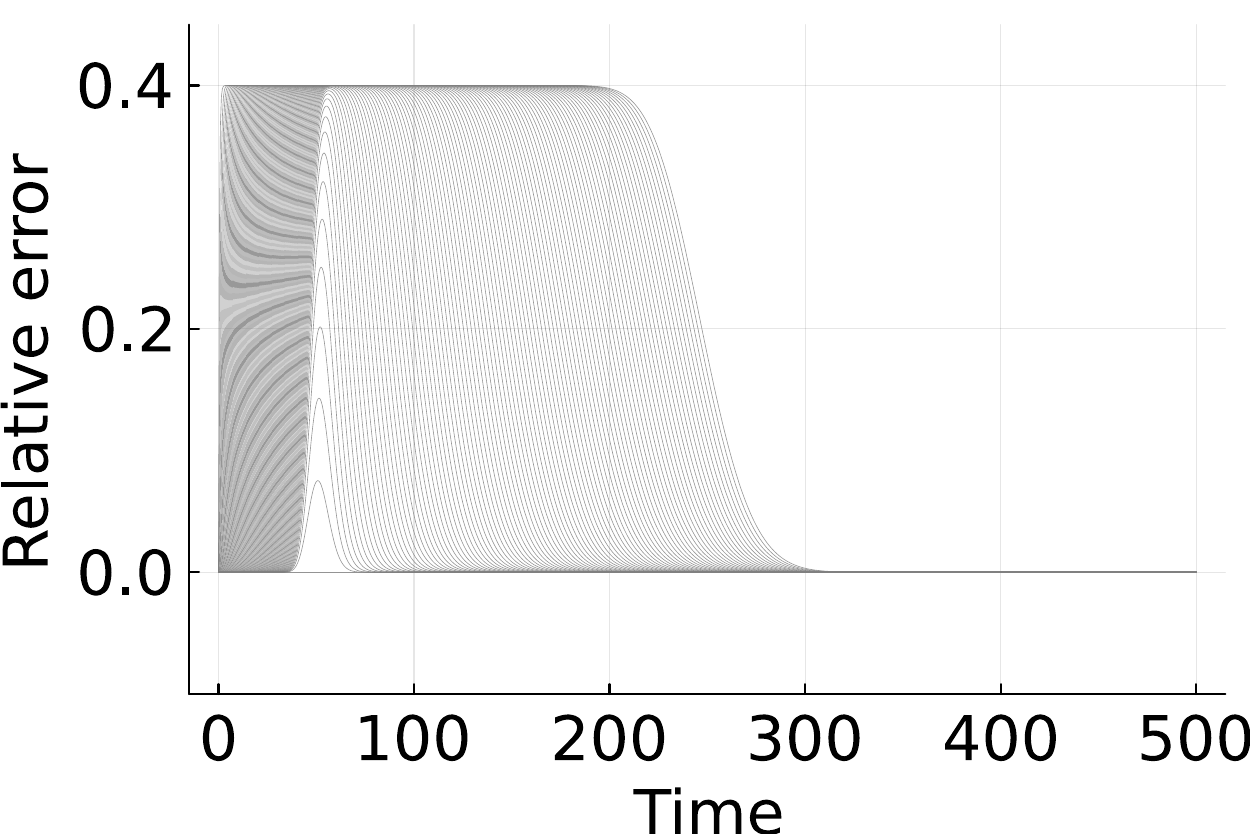}   
    \caption{Serial { consensus with bidirectional} feedback and $N=100$.}
    \label{fig:undir100_serial}
\end{subfigure}
    \caption{Simulation of the initial value response to $x=0$, $\dot{x}_{i\neq1}(0)=0$, and $\dot{x}_1(0)=1$. The conventional consensus system \eqref{eq:undir_conv} is considered with $r_1=2.5$ and $r_0=1$, while for the serial consensus system \eqref{eq:undir_serial}, $p_1=2$ and $p_2=0.5$ are used. The results illustrate that for some choices of graph Laplacians the serial and conventional consensus can have comparable performance.}
    \label{fig:bidirec_serial}
\end{figure}

\section{Conclusions and directions for future work}\label{sec:conclusions}

In this paper we have introduced the serial consensus controller which is a distributed formation controller that achieves scalable stability, performance, and robustness. Here, scalability refers to the fact that these properties are independent of the formation size. The performance result is particularly {important }%strong 
since it linearly bounds the \mbox{$\|\cdot\|_\infty$-gain} from the initial local errors and reference velocity deviation to the transient local errors and reference velocity deviations, measured in the same norm. {This quantity, rather than, for example an $L_2$ gain, is directly related to the control and performance objectives}. It is also worth noting that all of these results are achieved with only local relative measurements and linear feedback. The results are particularly interesting for large vehicle platoons where short inter-vehicle distances are desired and the transient behavior of the platoon is of great importance, {though there are strict topological constraints (typically, those of a directed string)}. But, by virtue of the generality of the presented results they could also be of interest for other networked systems, such as power grids, {sensor networks or multi-robot networks.} 

There are several interesting directions for future work. First, since the serial consensus may require an additional step of communication, {an interesting question is whether }  
this can be avoided, for instance through the use of local estimators. A second direction is to further investigate the robustness of the serial consensus, for instance {with respect} 
to time delays. Finally, {implementation of serial consensus on a physical system is an interesting next step.  }

\appendix
Here we will prove Theorem~\ref{thrm:scalable_stability} and Proposition~\ref{prop:q_implementable}. The results are restated for convenience.
\subsection{Proof of Theorem \ref{thrm:scalable_stability}}\label{app:thm1}
\begin{proof}
i) Any square matrix can be unitarily transformed to upper triangular form by the Schur traingularization theorem. Let $U_k L_kU_k^H= T_k$ be upper triangular. Then the block diagonal matrix $U=\mathrm{diag}(U_1, U_2, \dots U_{n})$ is a unitary matrix that upper triangularizes $A$ in \eqref{eq:serial_xiform}. For any triangular matrix the eigenvalues lie on the diagonal and this will be the eigenvalues of each $-L_k$. The result follows. 

 ii) First, 
    consider the closed-loop dynamics of \eqref{eq:serial_consensus} 
    $$X(s)=(sI+L_1)^{-1}(sI+L_2)^{-1}U_\mathrm{ref}(s).$$ 
    Since, $U_\mathrm{ref}$ is stable, we know that the limit $\lim_{s\rightarrow 0} U_\mathrm{ref}(s)=U_\mathrm{ref}(0)$ exists. To prove that the system achieves second-order consensus we want to show that 
    $$\lim_{t\rightarrow \infty}y(t)=\lim_{s\rightarrow 0}C(s)X(s)=0$$
    for some transfer matrix $C(s)$, which encodes the consensus states. But since the reference dependence is only related to $U_\mathrm{ref}(0)$, we can simplify the problem to only consider impulse responses, which %. But the impulse response 
    has the same transfer function as the initial value response where $\xi_2(0)=U_\mathrm{ref}(0)$. Therefore, WLOG, assume that $U_\mathrm{ref}(s)=0$ and an arbitrary initial condition $\mathbf{\xi}(0)=[\xi_1(0)^T,\xi_2^T(0)]^T.$ 
    The solution of~\eqref{eq:serial_xiform} is given by $\exp(At)\mathbf{\xi}(0)=S\exp(J(A)t)S^{-1}\mathbf{\xi}(0)$ where $J(A)$ is the Jordan normal form of $A$ and $S$ is an invertible matrix. From \ref{prop:high_i}) and the diagonal dominance of the graph Laplacians we know that all eigenvalues of $A$ lie in the left half plane. By  Assumption~\ref{ass:connected_tree} it follows that the zero eigenvalue for each~$L_k$ is simple. Now we prove that these two zero eigenvalues correspond to a Jordan block of size $2$. Let $\mathbf{e}_1=\bmat{\mathbf{1}^T& 0}^T$ and $\mathbf{e}_2=\bmat{0 &\mathbf{1}^T}^T$. Then $e_1$ is an eigenvector since $A\mathbf{e}_1=0$. Now, since $e_1=Ae_2$ combined with $e_1$ and $e_2$ being linearly independent it follows that they form a Jordan block of size $2$ with an invariant subspace spanned by the vectors $e_1$ and $e_2$. All other eigenvectors make up an asymptotically stable invariant subspace, it follows that $\xi(t)$ will converge towards a solution in $\mathrm{span}(e_1,e_2)$ and thus $\lim_{t\rightarrow \infty} \xi_k(t) =\alpha_k(t) \mathbf{1}$. From $x(t)=\xi_1(t)$ we get $\lim_{t\rightarrow \infty}x(t)=\alpha_1(t)\mathbf{1}$, and furthermore, since 
    $$\dot{x}=\dot{\xi}_1= -L_1\xi_1 +\xi_{2}\rightarrow \xi_{2}\text{ as } t\rightarrow \infty,$$
    it follows that $\lim_{t\rightarrow \infty}\dot{x}(t)= \alpha_{2}(t) \mathbf{1}$, which shows that the system achieves second-order consensus.
\end{proof}

\subsection{Proof of Proposition 2.}\label{app:prop}
\begin{proof}
    To prove this we must show that both $A_0,A_1\in\mathcal{A}^2(W,c')$  where $A_0=-L_2L_1$ and $A_1=-(L_2+L_1)$. First, $A_1$ and $A_2$ are shown to represent relative feedback. Since $L_1,L_2\in \mathcal{A}^1(W,c)$, it holds that $A_0\mathbf{1}=-L_2L_1\mathbf{1}=-L_20=0$ and similarly $A_1\mathbf{1}=-(L_2+L_1)\mathbf 1=0$. 
    
    Second, we show that the gain is bounded. For the positional feedback we have
    $$\|A_0\|_\infty=\|-L_2L_1\|_\infty\leq \|L_2\|_\infty\|L_1\|_\infty\leq c^2,$$
    which followed from the submultiplicativity of the induced norm. For the velocity feedback
    $$\|A_1\|_\infty=\|-(L_2+L_1)\|_\infty\leq \|L_2\|_\infty+\|L_1\|_\infty\leq 2c,$$
    where the triangle inequality was utilized. Let $c'=\max\{2c,c^2\}$, then clearly it holds true that both $\|A_0\|_\infty,\|A_1\|_\infty \leq c'$. 

    Finally, we consider the sparsity pattern. Since $W$ is non-negative, it follows that $W^2$ is also non-negative. Next, since adding non-negative elements to a matrix cannot remove any positive elements we get the following implications
    $$\left[I\!+\! W\!+\! W^2\right]_{i,j}\!\!=\!0\!\!\implies\!\!\!\left[I\!+\! W\right]_{i,j}\!\!=\!0\!\!\implies\!\!\![-L_1\!-\!L_2]_{i,j}\!\!=\!0,$$
    which follows from the definition of $L_1,L_2\in\mathcal{A}^1(W,c)$.
    
    To show that $A_0=-L_2L_1$ will be sparse, we note that any graph Laplacian can be written as $L=D-W$, where $D$ is a diagonal matrix and $W$ is non-negative. Thus we will consider the sparsity of
    $$-(D_2-W_2)(D_1-W_1)=-D_2D_1+D_2W_1+W_2D_1-W_2W_1.$$
    Since multiplication with a diagonal matrix preserves sparsity, it clearly holds true that the first three terms satisfy
    $$\left[I\!+\! W\right]_{i,j}\!\!=\!0\!\!\implies\!\!\![-D_2D_1+D_2W_1+W_2D_1]_{i,j}\!\!=\!0.$$
    For $W_2W_1$ we can introduce a non-negative matrix $E$ such that $[W_{1/2}+E]_{i,j}=0\iff W_{i,j}=0$ and in particular this construction ensures that $[(W_1+E)(W_2+E)]_{i,j}=0 \iff [W^2]_{i,j}$. Consider the following expanded product 
    $$(W_2+E)(W_1+E)=W_2W_1+W_2E+EW_1+E^2,$$
    where all terms are products of non-negative matrices and are therefore also non-negative. It thus holds that $[W^2]_{i,j}\implies [W_2W_1]_{i,j}=0.$ 
    Combining the results for $A_0$ shows that 
    $$\left[I\!+\! W\!+\!W^2\right]_{i,j}\!\!=\!0\implies \left[L_2L_1\right]_{i,j}=0.$$
    This concludes the proof.
\end{proof}

%%%

%\printbibliography 
\bibliographystyle{IEEETran}
\bibliography{references}
\end{document}